\documentclass[leqno,12pt]{article}
\usepackage{amssymb,amsfonts}
\usepackage{amsmath,latexsym}
\usepackage{amstext,epic,eepic,epsf}
\usepackage[dvips]{graphicx}

\newcommand{\bepr}{{\em Proof} } 
\newcommand{\enpr}{\hfill \rule{.5em}{.5em}}

\newcommand{\R}{{\mathbb R}}

\newcommand{\Z}{{\mathbb Z}}

\newtheorem{prop}{Proposition}[section] 
\newtheorem{thm}{Theorem}[section] 
\newtheorem{lemma}{Lemma}[section] 
\newtheorem{claim}{Claim}[section] 
\newtheorem{summ}{Summary}[section]

\begin{document}

\title{$L^2$-type Lyapunov functions for hyperbolic scalar conservation laws}

\author{Denis Serre \\ \'Ecole Normale Sup\'erieure de Lyon\thanks{U.M.P.A., UMR CNRS--ENSL \# 5669. 46 all\'ee d'Italie, 69364 Lyon cedex 07. France. {\tt denis.serre@ens-lyon.fr}}}

\date{\em Dedicated to Constantine Dafermos on the occasion \\ of his 80th birthday, with gratitude and admiration.}

\maketitle

\begin{abstract}
We prove the decay of the $L^2$-distance from the solution $u(t)$ of a hyperbolic scalar conservation law, to some convex, flow-invariant target sets. 
\end{abstract}

\paragraph{AMS classification}: 35L65, 35B35.

\paragraph{Key words}: Conservation laws, Lyapunov functions, shock waves.

\paragraph{Notations.} An $L^p$-norm is always denoted $\|\cdot\|_p$. The positive part of a real number $r$ is $r^+=\max\{r,0\}$. The differential of a $C^1$-function $\eta:\R^n\to\R$, computed at some point $a\in\R^n$, is ${\rm d}\eta(a)$.

\section{Introduction}

We are interested in decay/contraction properties for the flow defined by a scalar conservation law
\begin{equation}\label{eq:muld}
\partial_tu+{\rm div}_x\vec f(u)=0,\qquad x\in\R^d,\,t>0.
\end{equation}
The flux $\vec f:\R\to\R^d$ is a smooth function. Restricting to the natural notion of entropy weak solutions, Kru\v{z}kov's theory \cite{Kru} tells us that the Cauchy problem is well-posed in the class $L^\infty(\R^d)$. We therefore denote $(S_t)_{t\ge0}$ the semi-group defined by the flow of (\ref{eq:muld}). In one space dimension, the equation is written instead
\begin{equation}\label{eq:oned}
\partial_tu+\partial_xf(u)=0,\qquad x\in\R,\,t>0.
\end{equation}

When $b-a\in L^1(\R^d)$, we know that $S_tb-S_ta\in L^1(\R^d)$ as well, and $t\mapsto\|S_tb-S_ta\|_1$ is a non-increasing function. This contraction property does not extend to other $L^p$-norms, for $p>1$. We only have that if $a\in L^p(\R^d)$, then $t\mapsto\|S_ta\|_p$ is non-increasing, as a consequence of the fact that $s\mapsto|s|^p$ is a convex function. More generally, for every convex function $\eta:\R\to\R$, the admissible solution $u(t)=S_ta$ satisfies a differential ``entropy inequality''
$$\partial_t\eta(u)+{\rm div}_x\vec q(u)\le0,$$
where $\vec q$ is the entropy flux, defined by $\vec q\,'(s):=\eta'(s)\vec f'(s)$. If $\eta(0)=\eta'(0)=0$, this implies the monotonicity of
$$t\mapsto\int_{\R^d}\eta(u(t,x))\,dx,$$
whenever $\eta\circ a$ is integrable.

There are several reasons why we do not content ourselves with $L^1$-type results when studying conservation laws. One of them is that the ultimate goal of the theory is to consider {\em systems} of equations, such as the Euler system for inviscid compressible fluids, in several space dimensions. Then it is known (J. Rauch \cite{Rau}) that functional spaces such as $L^1(\R^d)$ and its avatar $BV(\R^d)$ are not appropriate: the Cauchy problem for linearized first-order hyperbolic systems is not well-posed in $L^p(\R^d)$ for $p\ne2$ (P. Brenner \cite{Bre}). This is why the most general existence and stability theorems have been established in classes build upon $L^2(\R^d)$, typically Sobolev spaces. For instance, it is known that the Cauchy problem for systems of conservation laws, endowed with a strongly convex entropy, is locally (in time) well-posed in $H^s(\R^d)$ whenever $s>1+\frac d2$\,, see \cite{Daf_b,BGS}. Notice that such solutions are continuously differentiable in space and time, thus exclude shock waves~; usually, they exist only for a finite time. Local existence and stability of shock waves, hence of discontinuous solutions, are proved by A. Majda in a space of piecewise $H^s(\R^d)$-functions, where the regularity parameter $s$ is even larger, see \cite{Maj1,Maj2}.

It is thus desirable to establish $L^2$-type results when discontinuous solutions are allowed. An early, notable step in this direction was done by C. Dafermos \cite{Daf} and R. J. DiPerna \cite{DiP}, who introduced the notion of {\em relative entropy} to prove weak-strong uniqueness and stability theorems. Roughly speaking, if a hyperbolic system of conservation laws is endowed with a strongly convex entropy $\eta$, and if the Cauchy problem admits a Lipschitz local-in-time solution, then the weak entropy solution is unique, equal to that one. The $L^2$-flavour of such a result is reminiscent to the fact that the relative entropy
$$\eta(u|v):=\eta(u)-\eta(v)-{\rm d}\eta(v)\cdot(u-v)$$
can be recast as a quadratic form $(u-v)^TA(u,v)(u-v)$ where the symmetric matrix $A$, which depends continuously upon $u$ and $v$, is positive definite. We notice that the uniqueness result is not associated with a contraction property, but follows from a Gronwall inequality. In particular, the time variable is present through an exponential factor.

When the reference solution $v$ is merely piecewise Lipschitz continuous, the relative entropy approach fails, because even a small disturbance in the initial data induces a change of the shock velocity and results in a rather large $L^2$-error. When $d=1$, this happens already when $v(t,x):=\phi(x-\sigma t)$ is a pure shock of (\ref{eq:oned}), defined by
\begin{equation}
\label{eq:pur}
\phi(y)=\left\{\begin{array}{lcr}
u_- & \hbox{if} & y<0, \\
u_+ & \hbox{if} & y>0.
\end{array}\right.
\end{equation}
This situation
motivated N. Leger and A. Vasseur \cite{Leg,LV} to enrich the approach by relative entropy, allowing the reference solution to be translated in space by a time-dependent shift. For instance Leger proved that if the flux $f$ is convex, then for every initial disturbance $a-\phi\in L^2(\R)$, the $L^2$-distance of $u(t)=S_ta$ to the set $\cal P$ of all translations of $\phi$,
$$d_2(u(t);{\cal P})=\inf_{h\in\R}\|u(t)-\phi(\cdot-h)\|_2,$$
is a non-increasing function of time.

The assumption that the flux be convex cannot be removed in Leger's theorem. When we turn towards multi-dimensional scalar conservation laws, this is a weakness, as we face directional fluxes $\xi\cdot\vec f$, for every unit vector $\xi$, which govern the propagation of planar waves $u(t,x)=U(t,\xi\cdot x)$. Unless the conservation law is essentially one-dimensional --~that is $\vec f(s)=f_0(s)\vec V$ for a fixed vector $\vec V$ --, there exist directions in which the flux is neither convex nor concave. Therefore there is no hope to extend Leger's result to the multi-D context. If we wish instead to drop the restriction that the flux be convex, we must pay a price, by changing $\cal P$ into a larger target set. A key remark is that $\cal P$, from a geometrical perspective, is not so nice: the curve $h\mapsto\phi_h:=\phi(\cdot-h)$  is not differentiable. Instead, it is of H\"older class with exponent $\frac12$\,:
$$\|\phi_k-\phi_h\|_2=|u_+-u_-|\cdot|h-k|^{1/2}.$$
This suggests to replace $\cal P$ by a smoother target set. Of course, we wish to keep its nice properties, namely the translation invariance -- because the semi-group itself is translation invariant, -- and its invariance under the flow. Our choice will be to replace $\cal P$ by its closed convex hull, which we describe below. It seems to us that $L^2$-stability properties of closed target set ${\cal Q}\subset L^\infty(\R^d)$ is related to the following properties
\begin{quote}
\begin{itemize}
\item convexity,
\item translation invariance~: if $a\in\cal Q$ and $h\in\R^d$, then $a_h:=a(\cdot-h)\in\cal Q$,
\item flow invariance~: if $a\in\cal Q$ and $t>0$, then $S_ta\in\cal Q$.
\end{itemize}
\end{quote}
Our main result below is that if $a\in \phi+L^2(\R)$, where $\phi$ is as in (\ref{eq:pur}), then the $L^2$-distance $d_2(S_ta;{\rm conv}({\cal P}))$ is a non-increasing function of time. We emphasize that {\em we do not assume} the convexity of the flux function $f$. Remarkably enough, we do not even need that the pure discontinuity $\phi$ be an admissible shock. Instead, remarking that the closed convex hull ${\cal M}:={\rm conv}({\cal P})$ coincides with the set of functions $a\in\phi+L^2(\R)$  that are monotonous from $u_-$ to $u_+$, our result completes the well-known fact that $S_t$ preserves the monotonicity.
\begin{thm}[$d=1$.]\label{th:mon}
Let $u_-,u_+\in\R$ be given. Let $\cal M$ be the set of monotonous functions $b$ over $\R$, such that $b-u_\pm\in L^2(\R^\pm)$ (in particular $b(\pm\infty)=u_\pm$).

If $a-u_\pm\in(L^1\cap L^\infty)(\R^\pm)$, then the $L^2$-distance $t\mapsto d_2(S_ta;{\cal M})$ is a non-increasing function.
\end{thm}

Of course, if the end states $u_\pm$ coincide, then ${\rm conv}({\cal P})$ is the singleton $\{a\equiv u_+\}$ and the result is nothing but the well-known decay of $t\mapsto\|u(t)-u_+\|_2$. Likewise, the theorem implies the invariance of $\cal M$ under the flow, a property which follows immediately from the comparison principle.

\bigskip

Our second example is multi-dimensional. Given $r>0$, our target set ${\cal Q}_r$ is the  intersection of $L^\infty(\R^d)$ with the $L^1$-ball of radius $r$. This too satisfies the three properties listed above. We notice that ${\cal Q}_r\subset L^2(\R^d)$.
\begin{thm}[$d\ge1$.]\label{th:ball}
Let ${\cal Q}_r$ be the $L^1$-ball of radius $r>0$ in $L^\infty(\R^d)$. 

If $a\in (L^2\cap L^\infty)(\R^d)$, then $t\mapsto d_2(S_ta;{\cal Q}_r)$ is non-increasing. 
\end{thm}

Both results are rather non-trivial. Even the projections upon either $\cal M$ or ${\cal Q}_r$ are not~; we describe them in detail in the next sections. That these theorems hold true reveals that we do not yet understand completely the nature of the semi-group $(S_t)_{t\ge0}$ for an arbitrary scalar conservation law.

\bigskip

As a remark, let us mention a few other situations that are more or less trivial:
\begin{itemize}
\item For $r>0$, the  intersection ${\cal B}_r$  of $L^\infty(\R^d)$ with the $L^2$-ball of radius $r$ satisfies the three properties listed above. Yet the distance $d_2(a;{\cal B}_r)$ equals $(\|a\|_2-r)^+$, thus the decay of $t\mapsto d_2(u(t);{\cal B}_r)$ is an obvious consequence of that of $\|u(t)\|_2$.
\item Let $I=[u_-,u_+]$ be a closed interval. The set ${\cal K}_I\subset L^\infty(\R^d)$ of functions which take values in $I$ satisfies the three properties listed above. The $L^2$-projection of $u$ over ${\cal K}_I$ is nothing but $\pi_I\circ u$, where $\pi_I:\R\to I$ is the usual projection. Thus
$$d_2(a;{\cal K}_I)=\left(\int_\R({\rm dist}(a(x);I))^2\,dx\right)^{\frac12}.$$
Since $\eta:=({\rm dist}(\cdot;I))^2$ is a convex function, the integral in the right-hand is non-increasing in time when applied to $u(t)=S_ta$, as a result of an entropy inequality.

Notice that if the flux displays enough non-linearity, and if $a\in \bar u+L^2(\R^d)$ for some $\bar u\in(u_-,u_+)$, we expect that $S_ta$ enters in ${\cal K}_I$ after some finite time. In the case of  the so-called multi-D Burgers equation, this property follows from the algebraic decay of $\|S_ta-\bar u\|_\infty$, proven by L. Silvestre and the author \cite{DSLS}.
\item  The replacement of the $L^2$-distance by the $L^1$-distance is just old stuff. If $\cal Q$ is  positively invariant, then the $L^1$-distance of $S_ta$ to ${\cal Q}$ is non-increasing in time, because $S_t$ is $L^1$-contracting.
\end{itemize}

\bigskip

\paragraph{Outline of the paper.} Section \ref{s:mainmon} is two-fold. On the one hand, it describes the effect of the projection over monotone functions. On the other hand, it displays the calculation behind Theorem \ref{th:mon} when the projection of the solution $u(t)$ behaves in a regular way in terms of the time variable. Section \ref{s:prfmon} is the technical part of the proof of Theorem \ref{th:mon}, where we succeed to reduce the analysis to the regular situation studied before. Surprisingly enough, a key argument pertains to Real Algebraic Geometry.  The proof of Theorem \ref{th:ball} is presented in Section \ref{s:prfball}. Because the projection over an $L^1$-ball is somehow a simpler operation, we can use the full strength of the kinetic formulation.

\bigskip

\paragraph{Acknowledgement.} I am indebted to Marie-Fran\c{c}oise Roy, who guided me in the realm of Real Algebraic Geometry.

\section{Main results towards Theorem \protect\ref{th:mon}}\label{s:mainmon}

Without loss of generality, we shall suppose $u_-<u_+$. The denote $\phi$ the pure discontinuity defined by (\ref{eq:pur}), even if it is not an admissible shock wave.

Recall that given a function $\psi:\R\to\R$, its {\em lower convex envelop} is the maximal convex function  $\rho\le\psi$. It is also the upper bound of the family of affine functions $\chi\le\psi$. If there does not exist such functions $\chi$, then we have $\rho\equiv-\infty$.

\subsection{Projection over $\cal M$}\label{ss:proj}

Let us denote $\pi$ the $L^2$-projection over $\cal M$, the set of functions $a\in\phi+L^2(\R)$ that are monotonous.
\begin{prop}\label{p:K}
Let $w\in(\phi+L^2(\R))\cap L^\infty(\R)$ be given. Then $\pi w=\rho'$ where $\rho$ is the lower convex envelop of a primitive $\psi$ of $w$.
\end{prop}

\bigskip

\bepr

By construction, $\rho'$ belongs to $\cal M$. Since this set is convex, and the projection is taken with respect to a Hilbertian norm, it suffices to prove that for every $g\in \cal M$, one has
\begin{equation}
\label{eq:topr}
\int_\R(w-\rho')(\rho'-g)\,dx\ge0.
\end{equation}
The open set $A=\{x;\rho(x)<\psi(x)\}$ is a union of disjoint intervals $(x_j,y_j)$. Mind that we do not exclude the possibility of a semi-infinite interval. Away from $A$, one has $\rho'=\psi'=w$. Thus the left-hand side of (\ref{eq:topr}) equals
$$\int_A(w-\rho')(\rho'-g)\,dx=\sum_j\int_{x_j}^{y_j}(\psi-\rho)'(\rho'-g)\,dx.$$
In a given integral of the right-hand side above, $\rho'$ is a constant $c_j$, the slope of the bi-tangent to the graph of $\psi$ between $x_j$ and $y_j$. Defining $h=g-c_j$ and $\theta=\psi-\rho$, we see that $\theta\ge0$ in $(x_j,y_j)$ is such that $\theta(y_j)=\theta(x_j)=0$, and $h$ is non-decreasing. Integrating by parts, we have
$$\int_{x_j}^{y_j}\theta'(x)h(x)\,dx=-\int_{x_j}^{y_j}\theta(x)h'(x)\,dx\le0.$$
Hence each term of the sum is $\ge0$. This proves (\ref{eq:topr}) and the statement.

\enpr

\bigskip

\paragraph{Remarks.} 
\begin{itemize}
\item Answering to a question raised on Mathoverflow\copyright, Willie Wong found the close formula for our projection
$$\pi w(x)=\inf_{z>x}\sup_{y\le x}\frac1{z-y}\,\int_y^zw(s)\,ds.$$
\item Let $\bar x$ be any of the points $x_j$ or $y_j$ in the proof above. Because $\rho$ is convex, $\rho\le\psi$ and $\rho(\bar x)=\psi(\bar x)$, we have $\psi'(\bar x-0)\le\rho'(\bar x-0)\le\rho'(\bar x+0)\le\psi'(\bar x+0)$. In other words the left and right limits of $w$ at $\bar x$ satisfy $w_\ell\le w_r$.
\end{itemize}

\subsection{The regular case: heuristic calculation}\label{ss:reg}

We now consider an admissible solution of (\ref{eq:oned}), with $u(t)\in\phi+L^2(\R)$. At each time, the primitive of the projection $\pi u(t)$ described in the previous paragraph differs from that of $u(t)$ over an open subset. We speak of the {\em regular} case when this open set depends smoothly upon $t$ on some time interval $(t_1,t_2)$. In particular, the limit points $x_j(t)$ and $y_j(t)$ are well-defined continuous and piecewise differentiable functions. In this situation, we prove a slightly more general result than just the decay of the $L^2$-distance to $\cal M$~:
\begin{prop}
\label{p:distM}
Let $a\in {\cal M}\cap BV(\R)$ be an initial data and $u(t):=S_ta$. Suppose that the projection onto $\cal M$ is regular for $t\in (t_1,t_2)$. Then, for every $C^2$-convex function $\eta$, the expression
$$\Delta(t)=\int_\R\eta(u(t)|\pi u(t))(x)\,dx,$$
is non-increasing in time.
\end{prop}

\bepr

Denote $(x_j(t),y_j(t))$ the disjoint intervals where the primitives of $u(t)$ and $\pi u(t)$ differ from each other. Then
$$\Delta(t)=\sum_j\int_{x_j}^{y_j}\eta(u(t,x)|c_j(t))\,dx,$$
where $\pi u(t)\equiv c_j(t)$ over $(x_j(t),y_j(t))$. 

We recall that $TV(u(t,\cdot))\le TV(a)$ and thus left and right traces $u(t,x\pm0)$ are well-defined. At every point $x_j$, we denote $v_{j,\ell/r}$ the left and right values $u(t,x_j\pm0)$. Likewise $w_{j,\ell/r}:=u(t,y_j\pm0)$. If there is no ambiguity about the point, then we write instead $u_{\ell/r}$. By the remark in the previous paragraph, we always have $u_\ell\le u_r$.

Let us differentiate $\Delta$ ~:
\begin{eqnarray*}
\dot\Delta & = & \sum_j\left( \dot y_j\eta(w_{j,\ell}|c_j)-\dot x_j\eta(v_{j,r}|c_j)+\int_{x_j}^{y_j}\partial_t(\eta(u|c_j))\,dx\right) \\
& \le & \sum_j\left(\dot y_j\eta(u_\ell|c_j)-\dot x_j\eta(u_r|c_j)-\int_{x_j}^{y_j}(\dot c_j\eta''(c_j)(u-c_j)+\partial_xq_j(u))\,dx\right)
\end{eqnarray*}
where $q_j$ is the entropy flux associated with the convex entropy $\eta(\cdot|c_j)$. We notice that the factor of $\dot c_j$ cancels\footnote{This is the reason why we choose a relative entropy, and not an arbitrary integrand $G(u(t),\pi u(t))$.} because in the projection we have
$$\int_{x_j}^{y_j}(u-c_j)\,dx=0.$$

Since $q_j'(s)=(\eta'(s)-\eta'(c_j))f'(s)$, we have
\begin{eqnarray*}
\dot\Delta &  \le & \sum_j\left(\dot y_j\eta(u_\ell|c_j)-\dot x_j\eta(u_r|c_j)+\int_{w_{j,\ell}}^{v_{j,r}}(\eta'(s)-\eta'(c_j))f'(s)\,ds\right) \\
& &   =\sum_j(A_{j}\dot x_j+B_{j}+C_j\dot y_j+D_j),
\end{eqnarray*}
where we denote
\begin{eqnarray*}
A_{j} = -\eta(v_{j,r}|c_j), & & B_{j} = \int_{c_j}^{v_{j,r}}(\eta'(s)-\eta'(c_j))f'(s)\,ds, \\
C_j = \eta(w_{j,\ell}|c_j), & & 
D_j = \int_{w_{j,\ell}}^{c_j}(\eta'(s)-\eta'(c_j))f'(s)\,ds.
\end{eqnarray*}

Let us look at the factor $A_{j}\dot x_j+B_{j}$. There are two cases, whether $u$ is continuous at $x_j$ or not. If it is not, then $u_\ell<c_j<u_r=v_{j,r}$, but then $u(t)$ displays a shock along $x=x_j$, so that Rankine--Hugoniot gives
$$\dot x_j=\frac{[f]}{[u]}\,.$$
Because the shock is increasing, we also have the entropy criterion that the graph of $f$ lies above its chord over $(u_\ell,u_r)$, denoted $s\mapsto {\rm Ch}(s)$. Since $\eta$ is convex, we have
\begin{eqnarray*}
B_{j} & = & (\eta'(u_r)-\eta'(c_j))f(u_r)-\int_{c_j}^{u_r}\eta''(s)f(s)\,ds \\
& \le & (\eta'(u_r)-\eta'(c_j))f(u_r)-\int_{c_j}^{u_r}\eta''(s){\rm Ch}(s)\,ds \\
& & = (\eta'(u_r)-\eta'(c_j))(\underbrace{f(u_r)-{\rm Ch}(u_r)}_{=0})+\int_{c_j}^{u_r}(\eta'(s)-\eta'(c_j)){\rm Ch}'(s)\,ds \\
& & =\frac{[f]}{[u]}\,\int_{c_j}^{u_r}(\eta'(s)-\eta'(c_j))\,ds  =\dot x_j\eta(u_r|c_j),
\end{eqnarray*}
because the slope of the chord is precisely the ratio $[f]/[u]$. We deduce that $A_{j}\dot x_j+B_{j}\le0$ in this case.

There remains the continuous case, where $u_\ell=u_r=c_j$. Here $A_{j}$ and $B_{j}$ vanish separately, so that again $A_{j}\dot x_j+B_{j}=0$.

The calculation is similar for the contribution $C_j\dot y_j+D_j$. We conclude that $\dot\Delta\le0$. In other words $t\longmapsto \Delta(t)$
is non-increasing.

\enpr

\section{Proof of Theorem \protect\ref{th:mon}}\label{s:prfmon}

We wish to apply Proposition \ref{p:distM} with the convex entropy $\eta(s)=s^2$. Because the solution $u(t)$ depends also of the initial data $a\in\phi+L^2(\R)$ and upon the flux $f$, we denote
$$\Delta_f(t;a)=\int_\R[(\pi u(t)-u(t))(x)]^2dx=(d_2(u(t);{\cal M}))^2.$$
Following the notation of Paragraph \ref{ss:proj}, we denote $A(t)$ the open set on which the primitive of $u(t,\cdot)$ differs from its lower convex envelop.

It is unclear whether the exact calculation of Paragraph \ref{ss:reg} can be applied directly to an entropy solution of (\ref{eq:oned}). It might happen that the structure of $A(t)$ varies so much that the calculation is not justified on any time interval. The combinatorial structure of the set of intervals could change drastically, infinitely many times.

\bigskip

We shall proceed as follows. By means of continuity arguments, prove that it suffices to consider $BV$ data and polynomial fluxes $f$~; see Summary \ref{sm:pol}. An other continuity argument allows us to limit ourselves to 
approximate solutions that are exact solutions between times steps, at which they are $L^2$-projected over piecewise constant functions~; see Claim \ref{c:app}. Such approximations can be obtained by the Godunov or Lax--Friedrichs schemes. We thus turn towards the case where the data, at some time step, is piecewise constant~; the corresponding solution concatenates Riemann problems, each one obeying an explicit formula involving an envelop. In particular $A(t)$ is always the union of finitely many disjoint intervals.
Because the flux is now a polynomial function, the solution of each Riemann problem can be expressed in terms of some semialgebraic set. Taking the convex envelop, as mentionned in Proposition \ref{p:K}, preserves this property. Then the Tarski--Seidenberg theorem ensures that the times at which a recombination occurs in the structure of $A(t)$, are finitely many. This allows us to apply piecewisely the  calculation of Paragraph \ref{ss:reg}.

\subsection{First reductions}

\paragraph{Dependence upon $a$.} Suppose that two data $a_1,a_2\in\phi+L^2(\R)$ take values in some bounded interval $[-M,M]$, and are such that $a_2-a_1\in L^1(\R)$. The corresponding solutions satisfy
$$\|u_j(t)\|_\infty\le M,\qquad\|u_2(t)-u_1(t)\|_1\le\|a_2-a_1\|_1,$$
which imply together
$$\|u_2(t)-u_1(t)\|_2^2\le2M\|a_2-a_1\|_1.$$
Since the distance to $\cal M$ is a $1$-Lipschitz function, we infer
$$\left|d_2(u_2(t);{\cal M})-d_2(u_1(t);{\cal M})\right|\le\sqrt{2M\|a_2-a_1\|_1\,}\,.$$
The functional $a\mapsto\Delta_f(t;a)$ is thus $L^1$-continuous over the set of data $a\in\phi+L^2(\R)$ with a prescribed pointwise bound $M$. 

Since the pointwise limit of non-increasing functions is non-increasing, it is enough to prove Theorem \ref{th:mon} for an $L^1$-dense subset of data. We shall therefore restrict our analysis to data $a\in BV(\R)$ which coincide with $\phi$ away from a bounded interval.

\paragraph{Time continuity.}
When $a\in BV(\R)\cap(\phi+L^2(\R))$, we know that 
$$\|u(t+h)-u(t)\|_1\le h\,{\rm Lip}\left(f|_{[-M,M]}\right)\,TV(a),$$
where $M=\|a\|_\infty$. Again, this gives $\|u(t+h)-u(t)\|_2\le C\sqrt h$ for some finite constant $C$, hence 
\begin{lemma}\label{l:dtwocont}
If $a\in BV(\R)\cap(\phi+L^2(\R))$, then $t\mapsto d_2(u(t);\cal M)$ is (H\"older) continuous.
\end{lemma}

\paragraph{Dependence upon the flux.}
We recall Lemma 11.1.1 of \cite{Ser_Cam2}~: Let $u$ be the entropy solution of (\ref{eq:oned}) with data $a\in BV(\R)$, and let $v$ be the entropy solutions of another equation $\partial_tv+\partial_xg(v)=0$, corresponding to the same initial data, then we have
$$\|v(T)-u(T)\|_1\le\int_0^TTV((g-f)\circ v(t))\,dt.$$
Combined with $TV(F\circ v)\le{\rm Lip}\left(F|_{[-M,M]}\right)\,TV(a)$, where $M=\|a\|_\infty$, this yields the ${\rm Lip}$-$L^1$ continuity of the map $f\mapsto u(t)$. With the same trick as above, we conclude that $f\mapsto\Delta_f(t;a)$ is continuous over ${\rm Lip}(-M,M)$, whenever $a\in BV(\R)\cap(\phi+L^2(\R))$ with $\|a\|_\infty\le M$. 

Thanks to this continuity property, we may restrict our study to fluxes that belong to a dense subspace of $C^1([-M,M])$. Applying the Stone--Weierstrass theorem (to $f'$ instead of $f$), we may restrict to {\em polynomial fluxes}.

\begin{summ}\label{sm:pol}
We only need to prove Theorem \ref{th:mon} when $a\in BV(\R)$ is such that $a-\phi$ is compactly supported, and the flux $f$ is a polynomial function.
\end{summ}

\paragraph{Approximate solutions.}
Our next remark is that the entropy solution of (\ref{eq:oned}) is the strong limit of the sequence of approximate solutions $u^h$ (with $h=\Delta x\to0+$ being the mesh size), generated by monotone difference schemes. Herebelow, we consider either the Godunov or the Lax--Friedrichs schemes, with a fixed CFL ratio 
\begin{equation}
\label{eq:CFL}
{\rm Lip}\left(f|_{[-M,M]}\right)\,\frac{\Delta t}h\,<\frac12\,.
\end{equation}
The convergence follows from Kuznetsov's estimate \cite{Kuz}~:
$$\|u^h(t)-u(t)\|_1\le C\sqrt{ht\,}\,TV(a).$$
Once again, we infer the $L^2$-convergence, whence
$$d_2(u(t);{\cal M})=\lim_{h\to 0+}d_2(u^h(t);{\cal M}),$$
so that
\begin{claim}\label{c:app}
To prove Theorem \ref{th:mon}, it suffices to verify that for every $h>0$, the function
$$t\mapsto d_2(u^h(t);{\cal M})$$
is non-increasing.
\end{claim}

\bigskip

Recall that the numerical scheme consists in alternating two operations. At each time step $t_k=k\Delta t$, $u^h(t_k-,\cdot)$ is interpolated by a piecewise constant function $u^h(t_k+,\cdot)$. This interpolation is nothing but the $L^2$-projection over the affine space of mesh-wise constant functions that tend to $u_\pm$ as $x\to\pm\infty$. When $k=0$, $u^h(0-)$ is simply the data $a$. Within an elementary time interval $(t_k,t_{k+1})$, $u^h$ is the (exact~!) entropy solution originating from the data $u^h(t_k+)$. It is obtained by concatenating solutions of Riemann Problems.

The interpolation step is the easy part of the analysis, as it does not involve the PDE at all:
\begin{lemma}\label{l:projdec}
Let $a\in\phi+L^2(\R)$ be given, and $\bar a$ be its $L^2$-projection over the  affine subspace of mesh-wise constant functions. Then 
\begin{equation}
\label{eq:dbara}
d_2(\bar a;{\cal M})\le d_2(a;{\cal M}).
\end{equation}
\end{lemma}

\bepr

Let us denote $I_j=((j-\frac12)h,(j+\frac12)h)$ the meshes. For definiteness, we consider the case of the Godunov scheme, where $j$ runs\footnote{If we worked with the Lax-Friedrichs scheme, $j$ would run over $\Z+\frac k2\,$, $k$ being the index of the time step.} over $\Z$. 


The projection, given by
$$(\bar w)|_{I_j}=\frac1h\,\int_{I_j}w(x)\,dx,$$
preserves the monotonicity: if $w\in\cal M$, then $\bar w\in\cal M$. Because this is an orthogonal projection, it is also a contraction. This implies 
$$\forall w\in{\cal M}, \qquad d_2(\bar u;{\cal M})\le\|\bar u-\bar w\|_2\le\|u-w\|_2.$$ 
Minimizing over $w\in\cal M$, we obtain (\ref{eq:dbara}).

\enpr

\subsection{Facts about Riemann problems}

The Riemann problem is the Cauchy problem for (\ref{eq:oned}) when the initial data is of the form
$$u(0,x)=\left\{\begin{array}{lcr}
v_-, & \hbox{if} & x<0, \\ v_+, & \hbox{if} & x>0,  
\end{array}\right.$$
where $v_\pm$ are two constants. The solution is self-similar, denoted
$$u(t,x)=R\left(\frac xt;v_-,v_+\right).$$

Because the initial data is monotonous, the solution is monotonous in the space variable as well. We may apply Matthias Kunic's formula \cite{Kun} (see also \cite{Ser_Cam1}, Proposition 2.5.1), which is dual to that of Lax~; it drops the assumption of a convex flux and asks instead for a monotonous data. We shall be concerned only by the non-decreasing case $v_-<v_+$, where the primitive $p(t,\cdot)$ of $u(t,\cdot)$ is given by
$$p(t,x)=\sup_s\inf_y\{s(x-y)-tf(s)+p_0(y)\}.$$
The primitive of $u(0,\cdot)$ being $p_0(y)=v_\pm y$ when $\pm y>0$, up to an additive constant, Kunik's formula yields
$$p(t,x)=tP\left(\frac xt\right),\qquad P(\xi):=\sup_{v_-\le s\le v_+}\{s\xi-f(s)\}.$$
Notice that the Legendre transform $P^*$ is the lower convex envelop of the restriction $f|_{[v_-,v_+]}$.

\subsection{End of the proof}

Because of Claim \ref{c:app} and Lemma \ref{l:projdec}, there remains to prove that $t\mapsto d_2(u^h(t);{\cal M})$ is non-increasing in each of the time intervals $(t_k+0,t_{k+1}-0)$. Translating in time, this amounts to prove Theorem \ref{th:mon} over $(0,\Delta t)$ whenever the initial data $a\in BV(\R)\cap (\phi+L^2(\R))$ is mesh-wise constant. Summary \ref{sm:pol} tells us that we may also assume on the one hand that only finitely many values
$$a_j:=a|_{I_j}$$
differ from $u_\pm$, and on the other hand the flux $f$ is a polynomial. We point out that the former constraint remains valid as time increases, the number of meshes where $a\ne\phi$ increasing only by $2$ at each time step. The main result of this paragraph is 
\begin{lemma}\label{l:main}
Denote $z$ the corresponding solution of the Cauchy problem associated with a mesh-wise constant initial data $\bar a=(a_j)_{j\in\Z}$.
Assume
\begin{itemize}
\item the flux $f$ is a polynomial,
\item $a_j\equiv u_-$ for $j<\!\!<-1$, while $a_j\equiv u_+$ for $j>\!\!>1$,
\item the CFL condition (\ref{eq:CFL}).
\end{itemize}
Then the interval $(0,\Delta t)$ splits into finitely many sub-intervals, in each of which the projection of $z(t)$ over $\cal M$ is regular, in the sense of Paragraph \ref{ss:reg}.
\end{lemma}

Proposition \ref{p:distM}, applied with $\eta(s)=s^2$, tells us that $d_2(u(t);{\cal M})$ is non-increasing within each of the sub-intervals mentionned in Lemma \ref{l:main}. Combining with the continuity stated in Lemma \ref{l:dtwocont}, we infer that it is non-increasing on the whole interval $(0,\Delta t)$. This ends the proof of Theorem \ref{th:mon}, provided we prove Lemma \ref{l:main}, which we do now.

\enpr

\bigskip

\bepr(of Lemma \ref{l:main}.)

The data $a$ is discontinuous at the grid points $x_{j+\frac12}$ for $j\in\Z$, which separate the states $a_j$ and $a_{j+1}$. Its primitive $p_0$ is continuous, piecewise linear. We denote $c_{j+\frac12}=p_0(x_{j+\frac12})$.

For every index $j\in\Z$, the solution $z$ in $(0,\Delta t)\times I_{j+\frac12}$ solves a Riemann Problem between the constant states $a_j$ and $a_{j+1}$. One has
$$z(t,x)=R\left(\frac{x-x_{j+\frac12}}t\,;a_j,a_{j+1}\right)=:Z_{j+\frac12}\left(\frac{x-x_{j+\frac12}}t\right).$$
Since this is a monotone function, the primitive $p(t,\cdot)$ of $z(t,\cdot)$ is either convex, or concave, on every mesh $I_{j+\frac12}$, depending on whether $a_j\le a_{j+1}$, or the opposite. It is given by
$$p(t,x)=:p_{j+\frac12}(t,x)=c_{j+\frac12}+tP_{j+\frac12}\left(\frac{x-x_{j+\frac12}}t\right)$$
where $P_{j+\frac12}$ is a primitive of $Z_{j+\frac12}$. Notice that $P_{j+\frac12}(\xi)=a_j\xi+d_j$ for $\xi<\!\!<0$ and $=a_{j+1}\xi+e_{j+1}$ for $\xi>\!\!>0$, where the integration constants satisfy $e_j=d_j$, because of the continuity of $p$ at $x_{j-\frac12}$.

Taking the lower convex envelop $q(t,\cdot)$ of $p(t,\cdot)$ is rather easy. Its graph differs from that of $p(t,\cdot)$ on bi-tangents, whose extremities belong to meshes where $p(t)$ is convex. Given two such meshes, there is at most one bi-tangent between them. In addition $p(t)$ needs to be concave somewhere in between, and thus the meshes may not be contiguous. Notice that a segment can be semi-infinite, meaning that it is tangent to the graph of $p(t)$ at a finite point and at $\pm\infty$. Since $z(t)\equiv\phi$ away from a compact interval, $q(t,\cdot)$ differs from $p(t,\cdot)$ on finitely many segments only.

Define $J\subset\Z$ the finite set of indices such that $a_j<a_{j+1}$. To determine the bi-tangents, we begin by selecting $i,j$ in $J$ such that $j-i\ge2$. Because $p(t)$ is convex in both $I_{i+\frac12}$ an $I_{j+\frac12}$, there is at most one bi-tangent whose tangency points belong to both meshes. Its slope being $\theta$, the tangency in $I_{j+\frac12}$ tells us that its equation is
$$q=\theta x-(p_{j+\frac12})^*(\theta)=\theta(x-x_{j+\frac12})+c_{j+\frac12}-t(P_{j+\frac12})^*(\theta).$$
Expressing the tangency in $I_{i+\frac12}$, we obtain an alternate equation of the bi-tangent:
$$q=\theta(x-x_{i+\frac12})+c_{i+\frac12}-t(P_{i+\frac12})^*(\theta).$$
Eliminating, we find that the slope $\theta$ is determined by the equation
\begin{equation}
\label{eq:biT}
(j-i)h\theta+t(P_{j+\frac12})^*(\theta)-t(P_{i+\frac12})^*(\theta)=c_{j+\frac12}-c_{i+\frac12}.
\end{equation}

Let us recall that $(P_{j+\frac12})^*$ is the lower convex envelop of the restriction ot the flux $f$ to $(a_j,a_{j+1})$. It is therefore a $C^1$-function on this interval, whose derivatives are derivatives of $f$. The existence of a bi-tangent necessitates that $(a_i,a_{i+1})\cap(a_j,a_{j+1})\ne\emptyset$. Because of the CFL condition, the left-hand side of (\ref{eq:biT}) is a uniformly increasing function of $\theta$. Thanks to the Implicit Function Theorem, the slope $\theta=\theta(t)$ of the bi-tangent is a $C^1$ function of time.

Our solution $z$ is thus {\em regular} on every time interval on which the bi-tangents depend continuously upon $t$. A recombination of $A(t)$ may occur only if two consecutive bi-tangents, corresponding to pairs $(i,j)$ and $(j,k)$, happen to coincide. Thus we are lead to study the occurences of tri-tangents to the graph of $p(t)$.

For a tri-tangent to occur at some $t\in(0,\Delta t)$, one needs a triple $i<j<k$ of elements of $J$. Then the slope $\theta$ satisfies
\begin{eqnarray}
\label{eq:trij}
(j-i)h\theta+t(P_{j+\frac12})^*(\theta)-t(P_{i+\frac12})^*(\theta) & = & c_{j+\frac12}-c_{i+\frac12}, \\
\label{eq:trjk}
(k-j)h\theta+t(P_{k+\frac12})^*(\theta)-t(P_{j+\frac12})^*(\theta) & = & c_{k+\frac12}-c_{j+\frac12}.
\end{eqnarray}
We may express the solutions of (\ref{eq:trij}) (respectively of (\ref{eq:trjk})) by $\theta=\Theta_{ij}(t)$ (resp. $\theta=\Theta_{jk}(t)$) where the functions $\Theta_{\cdot\cdot}$ are $C^1$. Thus $A(t)$ may recombine only at times such that $\Theta_{ij}(t)=\Theta_{jk}(t)$.
For a general flux, the solution set of this equation can be extremely complicated. But after our reductions, we need only to consider the case of a polynomial flux. 

Recall that $(P_{j+\frac12})^*$ is, up to an additive constant (the constants $d_j$ above), the lower convex envelop of the restriction $f_{j+\frac12}$ of $f$ to $(a_j,a_{j+1})$. Its calculation requires the computation of the bi-tangents to $f_{j+\frac12}$. This is an elimination in a system of algebraic equations. Selecting the relevant bi-tangents (those which are below the graph of $f_{j+\frac12}$) requires adding algebraic inequalities. The result of such operations is that the graph of $(P_{j+\frac12})^*$ is a semialgebraic set, meaning that it is a finite union of real sets defined by polynomial identities and polynomial inequalities. Since (\ref{eq:trij},\ref{eq:trjk}) is a polynomial system in $(\theta,t,(P_{i+\frac12})^*(\theta),(P_{j+\frac12})^*(\theta),(P_{k+\frac12})^*(\theta))$, its solutions $(\theta,t)$ form a semialgebraic set ${\cal S\!A}_{ijk}$. 

By the Tarski--Seidenberg Principle, the projection ${\cal S\!T}_{ijk}$ of ${\cal S\!A}_{ijk}$ on the time axis is still a semialgebraic set, see \cite{BCR} Chapter 5. It is therefore a finite union of points and intervals\footnote{These assertions tell us that semialgebraic sets form an {\em o-minimal} structure.}. This projection is precisely the set of times $t\in(0,\Delta t)$ at which a tri-tangent occurs, with tangencies in the meshes of indices $i,j,k$. If it occurs an isolated point, fine~! If instead it occurs along a time interval $(t_-,t_+)$, then we may ignore the tangency in the intermediate mesh $I_{j+\frac12}$, and consider this tri-tangent as a regular bi-tangent between $I_{i+\frac12}$ and $I_{k+\frac12}$. Eventually, replacing each of the segments of ${\cal S\!T}_{ijk}$ by its extremities, this set becomes equivalent, from the point of view of the regularity of $z$, to a finite set $\widetilde{\cal S\!T}_{ijk}$. 

Since the admissible triples $(i,j,k)$ are finitely many (because $J$ is finite), the union $\widetilde{\cal S\!T}$ of the sets $\widetilde{\cal S\!T}_{ijk}$ is still finite. It splits $(0,\Delta t)$ into finitely many sub-intervals, on which our solution is regular in the sense of Paragraph \ref{ss:reg}.

\enpr

\bigskip

This ends the proof of Theorem \ref{th:mon}.

\section{Proof of Theorem \protect\ref{th:ball}}\label{s:prfball}

The situation is now multi-dimensional. We consider data in $(L^2\cap L^\infty)(\R^d)$, a domain invariant under the action of the semi-group. The target set is the ball ${\cal Q}_r$ defined by
$$\|a\|_1:=\int_{\R^d}|a(x)|\,dx\le r.$$

\subsection{Projection over ${\cal Q}_r$}

If $r>0$, the intersection $D_r$ of the closed ball $B(0;r)$ in $L^1(\R^d)$, with $L^2(\R^d)$ is a closed convex subset of the latter. Let $\pi_r:L^2(\R^d)\rightarrow D_r$ be the projection according to the natural distance $d(v,w)=\|w-v\|_2$. The $L^2$-projection from $L^2\cap L^\infty$ onto ${\cal Q}_r$ is nothing but the restriction of $\pi_r$.

\begin{prop}\label{p:projnosgn}
If $v\in L^2(\R^d)$, then either $\pi_rv=v$ (if $\|v\|_1\le r$), or $\pi_rv=({\rm sgn}\,v)(|v|-s)^+$ where $s\ge0$ is determined by
$$\int_{\R^d}(|v|-s)^+dx=r,$$
if instead $\|v\|_1>r.$
\end{prop}

\bepr

We may assume the latter situation. If $s>0$, then
$$\int_{\R^d}(|v|-s)^+dx\le\int_{\{|v|>s\}}|v(x)|\,dx\le\frac1s\,\int_{\R^d} v(x)^2dx<\infty.$$
The map $I:s\mapsto\int(|v|-s)^+dx$ is non-increasing, ranging from $I(0+)=\|v\|_1$ (possibly infinite) to $I(\|v\|_\infty)=0$. It is actually Lipschitz continuous away from the origin, because is $0<s<t$, then
$$I(s)-I(t)\le(t-s){\rm meas}\{|v|>t\}\le\frac{\|v\|_2^2}{t^2}\,(t-s).$$
At last, $I$ is strictly monotonous over $(0,\|v\|_\infty)$ because if $0<s<t< |v|$, then $(|v|-t)^+-(|v|-s)^+=s-t<0$. There exists therefore a unique $s\in(0,\|v\|_\infty)$ such that $I(s)=r.$

Proving that $g:=({\rm sgn}\,v)(|v|-s)^+$ is the projection amounts to verifying $\langle v-g,g-h\rangle\ge0$ for every $h\in D_r$. But this quantity equals 
\begin{eqnarray*}
\int_{\{|v|>s\}}s(|g|-h{\rm sgn}\,v)\,dx-\int_{\{|v|\le s\}}vh\,dx & = & s\|g\|_1-s\int_{\{|v|>s\}}h\,{\rm sgn}\,v\,dx-\int_{\{|v|\le s\}}vh\,dx \\ &  \ge & sr-s\|h\|_1\ge0.
\end{eqnarray*}

\enpr

\subsection{Proof by the kinetic formulation}

We have the slightly more general result:
\begin{thm}\label{th:mdsgngen}
Let $a\in L^2(\R^d)$ be an initial data and $u(t):=S_ta$.
Let $\eta$ be a smooth even convex function with $\eta(0)=0$. Denote $v(t)=\pi_ru(t)$.
Then the expression
$$\Delta(t):=\int_{\R^d}\eta(u(t,x)-v(t,x))\,dx$$ is a non-increasing function of time.
\end{thm}

We emphasize the fact that the statement is valid only if $\eta$ is even (see the proof below), and that it does not involve a relative entropy. We notice also that, because $t\mapsto\|u(t)\|_1$ is non-increasing, the projection $\pi_r$ acts non-trivially for $t$ in some time interval $(0,T)$ (with $T$ possibly infinite or null), and then trivially for $t\ge T$, in which case $\Delta(t)=0$.

\bigskip

\bepr

We use the kinetic formulation of (\ref{eq:muld}). To this end, we recall the definition of the chi-function:
$$\chi(\xi;u)=\left\{\begin{array}{lcr}
1 & \hbox{if} & 0<\xi<u, \\
-1 & \hbox{if} & u<\xi<0, \\
0 & & \hbox{otherwise.}
\end{array}\right.$$

The function $u(t,x)$ is an entropy solution of (\ref{eq:muld}) if and only if there exists a non-negative bounded measure $m\in C(\R_\xi;{\cal M}(\R_+\times\R^d))$, such that the kinetic density $h(t,x,\xi):=\chi(\xi;u(t,x))$ satisfies the transport equation (see \cite{LPT}, or Theorem 3.2.1 of \cite{Per})
\begin{equation}
\label{eq:kinfor}
\partial_th+f'(\xi)\cdot\nabla_xh=\frac{\partial m}{\partial\xi}\,.
\end{equation}

We may restrict to the time interval $(0,T)$. We have
\begin{eqnarray*}
r & \equiv & \int_{\R^d}dx\left(\int_s^{+\infty}-\int_{-\infty}^{-s}\right)h(t,x,\xi)\,d\xi, \\
\Delta(t) & = & \int_{\R^d}\eta(\min(|u|,s)\,dx=\int_{\R^d}dx\int_{-s}^s\eta'(\xi)h(t,x,\xi)\,d\xi.
\end{eqnarray*}
Differentiating, there comes
\begin{eqnarray*}
\dot s\int_{\R^d}(h(t,x,s)-h(t,x,-s))\,dx & = & \int_{\R^d}dx\left(\int_s^{+\infty}-\int_{-\infty}^{-s}\right)\partial_th(t,x,\xi)\,d\xi, \\
\dot\Delta(t) & = & \int_{\R^d}dx\int_{-s}^s\eta'(\xi)\partial_th(t,x,\xi)\,d\xi \\
& & +\dot s\int_{\R^d}(\eta'(s)h(t,x,s)+\eta'(-s)h(t,x,-s))\,dx.
\end{eqnarray*}
Since $\eta$ is even, and thus $\eta'$ is odd, we can eliminate $\dot s$, to obtain
$$\dot\Delta(t) = \int_{\R^d}dx\int_{-s}^s\eta'(\xi)\partial_th(t,x,\xi)\,d\xi +\eta'(s) \int_{\R^d}dx\left(\int_s^{+\infty}-\int_{-\infty}^{-s}\right)\partial_th(t,x,\xi)\,d\xi.$$
Replacing $\partial_th$ by $\partial_\xi m-f'(\xi)\cdot\nabla_x h$ in the identity above, and then integrating by parts in the space variable, there remains
$$\dot\Delta(t) = \int_{\R^d}dx\int_{-s}^s\eta'(\xi)\partial_\xi m\,d\xi +\eta'(s) \int_{\R^d}dx\left(\int_s^{+\infty}-\int_{-\infty}^{-s}\right)\partial_\xi m\,d\xi.$$
The last term above is non-positive because on the one hand $\eta'(s)\ge0$ and on the other hand $m$ is a non-negative finite measure in $\xi$. hence
$$\dot\Delta(t) \le \int_{\R^d}dx\int_{-s}^s\eta'(\xi)\partial_\xi m\,d\xi  = -\int_{\R^d}dx\int_{-s}^s\eta''(\xi)m\le0.$$

\enpr

\bigskip

Remark that, integrating the latter in time, we obtain an estimate
$$\int_0^{+\infty}\int_{\R^d}\int_0^{s(t)}\eta''(\xi)m\le\,\Delta(0).$$
Choosing $\eta(u)=u^2$, letting $r\rightarrow0+$, which yields $s\rightarrow\|u(t)\|_\infty$, we recover the well-known inequality (see Proposition 3.2.3 of \cite{Per}))
$$\int_0^{+\infty}\int_{\R^d}\int_0^{+\infty}m\le\frac12\,\|a\|_2^2.$$


\begin{thebibliography}{00}





\bibitem{BCR} J. Bochnak, M. Coste, M.-F. Coste-Roy. {\em Real Algebraic Geometry}. Ergebnisse der Math. und ihrer Grenzgebiete {\bf36}, Springer-Verlag (1998).








\bibitem{BGS} S. Benzoni-Gavage, D. Serre. {\em Multi-dimensional hyperbolic partial differential equations. First-order systems and applications.} Oxford Math. Monographs. Oxford Univ. Press, Oxford (2007).




\bibitem{Bre} P. Brenner. The Cauchy problem for symmetric hyperbolic systems in $L^p$. {\em Math. Scand.}, {\bf19} (1966), pp 27--37.









\bibitem{Daf} C. Dafermos. The second law of thermodynamics and stability. {\em Arch. Rat. Mech. Anal.}, {\bf 70} (1979), pp 167--179.


\bibitem{Daf_b} C. Dafermos. {\em Hyperbolic conservation laws in continuum physics.} Grundlehren der mathematischen Wissenschaften {\bf325}. Springer-Verlag, Berlin (2000).


        


\bibitem{DiP} R. DiPerna. Uniqueness of solutions to hyperbolic conservation laws. {\em Indiana Univ. Math. J.}, {\bf 28} (1979), pp 137--188.
















\bibitem{Kun} M. Kunik. A solution formula for a non-convex scalar hyperbolic conservation law with monotone initial data. {\em Math. Methods Appl. Sciences}, {\bf16} (1993), pp 895--902.



        
        
\bibitem{Kru} S. Kru\v{z}kov. First order quasilinear equations with several independent variables (in Russian). {\em Mat. Sbornik (N.S.)}, {\bf 81 (123)} (1970), pp 228--255.
%


        
\bibitem{Kuz} S. Kuznetsov. Accuracy of some approximate methods for computing the weak solutions of a first-order quasi-linear equation. {\em USSR Comp. Math. and Math. Phys.}, {\bf16} (1976), pp 105--119.


\bibitem{Leg} N. Leger. $L^2$-stability estimates for shock solutions of scalar conservation laws using the relative entropy method. {\em Arch. Rat. Mech. Anal.}, {\bf 199} (2011), pp 761--778.
        
\bibitem{LV} N. Leger, A. Vasseur. Relative entropy and the stability of shocks and contact discontinuities for systems of conservation laws with non-BV perturbations. {\em Arch. Rat. Mech. Anal.}, {\bf 201} (2011), pp 271--302.
        


\bibitem{LPT} P.-L. Lions, B. Perthame, E. Tadmor. A kinetic formulation of multidimensional scalar conservation laws and related equations. {\em J. Amer. Math. Soc.}, {\bf 7} (1994), pp 169--191.


\bibitem{Maj1} A. Majda. {\em The stability of multidimensional shock fronts.} Mem. Amer. Math. Soc. {\bf41} (1983).


\bibitem{Maj2} A. Majda. {\em The existence of multidimensional shock fronts.} Mem. Amer. Math. Soc. {\bf43} (1983).
        

        





        


        


\bibitem{Per} B. Perthame.  {\em  Kinetic formulation of conservation laws}, Oxford lecture series in Math. \& its Appl. {\bf21}. Oxford (2002).

\bibitem{Rau} J. Rauch. BV estimates fail for most quasilinear hyperbolic systems in dimension greater than one. {\em Comm. Math. Phys.}, {\bf106} (1986), pp 481--484.  

        

        

        

        




\bibitem{Ser_Cam1} D. Serre. {\em Systems of conservation laws 1: Hyperbolicity, entropies, shock waves.} Cambridge Univ. Press (1999), Cambridge, UK.





\bibitem{Ser_Cam2} D. Serre. {\em Systems of conservation laws 2: Geometric structures, oscillations and initial boundary-value problems.} Cambridge Univ. Press (2000), Cambridge, UK.








      


\bibitem{DSLS} D. Serre, L. Silvestre. Multi-dimensional scalar conservation laws with unbounded initial data: well-posedness and dispersive estimates. {\em Arch. Rat. Mech. Anal.}, {\bf234} (2019), pp 1391--1411.

        

        



        

        

        



 

\end{thebibliography}
\end{document}